\documentclass[12pt,twoside]{amsart}
\usepackage{graphicx}
\usepackage{color}
\usepackage{amssymb}
\usepackage{amsmath}
\usepackage{amsthm}
\usepackage{amsfonts}
\usepackage{amsmath, amsfonts, amssymb}
\newtheorem{theorem}{Theorem}[section]

\newtheorem{definition}[theorem]{Definition}
\numberwithin{equation}{section}

\begin{document}

	\title{A Description of Automorphisms and Derivations of All Two-Dimensional Algebras over an Arbitrary Field} 
	\author{ Eshmirzayev Sh.$^1$, Bekbaev U.$^2$} 
	\thanks{{\scriptsize emails: shoxjahoneshmirzayev.95@mail.ru, uralbekbaev@gmail.com}}
	\begin{center} \address{$^1$  Tashkent Information Technology University, Tashkent, Uzbekistan;\\
			\smallskip 
			$^2$ Turin Polytechnic University in Tashkent, Tashkent, Uzbekistan}
	\end{center}
	\begin{abstract} A description of automorphisms and derivations of all two-dimensional algebras, considered up to isomorphism, over an arbitrary field is provided. In other words, these two problems are completely resolved in this paper.  	\end{abstract}
	\maketitle
	
	\section{Introduction} Description of automorphisms and derivation algebras of a given algebra $\mathbb{A}$ is considered an important part of the study of that algebra. In \cite{ABR}, such a description for all two-dimensional algebras over any field $\mathbb{F}$, considered up to isomorphism and under the assumption that every second- and third-degree polynomial over the base field has a root in it, is given. The main aim of this paper is to describe the automorphisms and derivations of such algebras, considered up to isomorphism, without any restrictions on the base field $\mathbb{F}$. In other words, we completely resolve the problem of describing automorphisms and derivations of two-dimensional algebras.
	
	Moreover, in contrast to \cite{ABR}, we use, whenever possible, certain relative invariants that significantly simplify the proof of the main result.
	
	\section{Preliminaries}
	
	Let $\mathbb{F}$ be any field. For matrices $A=(a_{ij})$ and $B$ over $\mathbb{F}$, as usual, $A\otimes B$ denotes the block matrix with blocks $(a_{ij}B)$.
	
	\begin{definition}
		A vector space $\mathbb{A}$ over $\mathbb{F}$ equipped with a multiplication $\cdot :\mathbb{A}\times \mathbb{A}\rightarrow \mathbb{A}$,	
		given by $(\mathbf{u},\mathbf{v})\mapsto \mathbf{u}\cdot \mathbf{v}$ such that
		\[(\alpha\mathbf{u}+\beta\mathbf{v})\cdot \mathbf{w}=\alpha(\mathbf{u}\cdot \mathbf{w})+\beta(\mathbf{v}\cdot \mathbf{w}),\ \ \mathbf{w}\cdot (\alpha\mathbf{u}+\beta\mathbf{v})=\alpha(\mathbf{w}\cdot \mathbf{u})+\beta(\mathbf{w}\cdot \mathbf{v})\] whenever $\mathbf{u}, \mathbf{v}, \mathbf{w}\in \mathbb{A}$ and $\alpha, \beta\in \mathbb{F}$, is called an algebra.
	\end{definition}
	
	\begin{definition}
		Two algebras $\mathbb{A}$ and $\mathbb{B}$ are called isomorphic if there exists an invertible linear map $\mathbf{f}:\mathbb{A}\rightarrow \mathbb{B}$ such that
		\[\mathbf{f}(\mathbf{u}\cdot_{\mathbb{A}} \mathbf{v})=\mathbf{f}(\mathbf{u})\cdot_{\mathbb{B}} \mathbf{f}(\mathbf{v})\]
		whenever $\mathbf{u}, \mathbf{v}\in \mathbb{A}$.
	\end{definition}
	
	\begin{definition}
		A linear map $\mathbf{g}:\mathbb{A}\rightarrow \mathbb{A}$ (respectively, $\mathbf{D}:\mathbb{A}\rightarrow \mathbb{A}$) is called an automorphism(derivation) if $\mathbf{g}$ is invertible and
		\[\mathbf{g}(\mathbf{u}\cdot \mathbf{v})=\mathbf{g}(\mathbf{u})\cdot \mathbf{g}(\mathbf{v}),\ (\mbox{respectively},\ \mathbf{D}(\mathbf{u}\cdot \mathbf{v})=\mathbf{D}(\mathbf{u})\cdot \mathbf{v}+\mathbf{u}\cdot\mathbf{D}(\mathbf{v})) \] 	whenever $\mathbf{u}, \mathbf{v}\in \mathbb{A}$.
	\end{definition}
	
	Let $\mathbb{A}$ be an $m$-dimensional algebra over $\mathbb{F}$, and let $e=(e_1,e_2,\ldots,e_m)$ be a basis of $\mathbb{A}$. Then the bilinear map $\cdot$ is represented by a matrix $A=(A^k_{ij})\in M(m\times m^2;\mathbb{F})$ as follows. Let $\mathbb{A}$ be $m$-dimensional algebra over $\mathbb{F}$ and $e=(e_1,e_2,...,e_m)$ its basis. Then the bilinear map $\cdot$ is represented by a matrix $A=(A^k_{ij})\in M(m\times m^2;\mathbb{F})$ as follows:
	\[\mathbf{u}\cdot \mathbf{v}=eA(u\otimes v),\] for $\mathbf{u}=eu,\mathbf{v}=ev,$
	where $e_ie_j=\sum_{k=1}^mA^k_{ij}e_k$, $u,v\in \mathbb{F}^m$ are column coordinate vectors of $\mathbf{u}$ and $\mathbf{v},$ respectively. 
	
	The matrix $A\in M(m\times m^2;\mathbb{F})$ defined above is called the matrix of structural constants (MSC) of $\mathbb{A}$ with respect to the basis $e$. Throughout the paper, we assume that a basis $e$ is fixed, and we do not distinguish between the algebra $\mathbb{A}$ (respectively, the $\mathbb{F}$-linear map $\mathbf{g}:\mathbb{A}\rightarrow \mathbb{A}$) and its MSC $A$ (respectively, the matrix $g\in Mat(m;\mathbb{F})$) with respect to $e$.
	
	Note that a linear map $\mathbf{g}:\mathbb{A}\rightarrow \mathbb{A}$ (respectively, $\mathbf{D}:\mathbb{A}\rightarrow \mathbb{A}$) is an automorphism (respectively, a derivation) if and only if $g\in GL(m;\mathbb{F})$(respectively, $D\in Mat(m;\mathbb{F})$) and the equality  \begin{equation}\label{1}gA=Ag^{\otimes 2},\ (\mbox{respectively},\ DA=A(D\otimes I+I\otimes D)\end{equation} holds, where $g$ (respectively, $D$) is the matrix of $\mathbf{g}$ (respectively, $\mathbf{D}$) with respect to the basis $e$, $I$ denotes the identity matrix of order $m$.
	
	Note that the isomorphism of two $m$-dimensional algebras, in terms of their MSCs, is equivalent to the following definition.
	
	\begin{definition}
		Two $m$-dimensional algebras $\mathbb{A}$ and $\mathbb{B}$ over $\mathbb{F}$, given by their matrices of structural constants $A$ and $B$, are said to be isomorphic if
		$B=gA(g^{-1})^{\otimes 2}$
		for some $g\in GL(m;\mathbb{F})$.
	\end{definition}
	
	From now on, we consider only the case $m=2$. For simplicity, we write
	\[A=\left(\begin{array}{cccc} \alpha_1 & \alpha_2 & \alpha_3 &\alpha_4\\ \beta_1 & \beta_2 & \beta_3 &\beta_4\end{array}\right)\]
	for an MSC, where $\alpha_1, \alpha_2, \alpha_3, \alpha_4, \beta_1, \beta_2, \beta_3, \beta_4$ are arbitrary elements of $\mathbb{F}$.
	
	Our proof of the main result relies on Classification Theorems 1, 2, and 3 proved in \cite{B}, together with the corrections presented in \cite{B1}. For the convenience of the reader, we present the corrected results below in an unified form.
	
	\begin{theorem}Any non-trivial two-dimensional algebra over a field $\mathbb{F}$ with $\mathrm{Char}(\mathbb{F})\neq 2,3$ is isomorphic to exactly one of the following algebras, listed by their matrices of structural constants:\begin{itemize}
			\item
			$ A_{1}(\mathbf{c})=\left(
			\begin{array}{cccc}
				\alpha_1 & \alpha_2 &1+\alpha_2 & \alpha_4 \\
				\beta_1 & -\alpha_1 & 1-\alpha_1 & -\alpha_2
			\end{array}\right),\ \mbox{where}\ \mathbf{c}=(\alpha_1, \alpha_2, \alpha_4, \beta_1) \in \mathbb{F}^4,$\\ 
			\item	$ A_{2}(\mathbf{c})=\left(
			\begin{array}{cccc}
				\alpha_1 & 0 & 0 & \alpha_4 \\
				1& \beta_2 & 1-\alpha_1&0
			\end{array}\right),\ \mbox{where}\ \mathbf{c}=(\alpha_1,\alpha_4, \beta_2)\in \mathbb{F}^3,\ \alpha_4\in \mathbb{F}^*,$ 
			\item	$ A_{3}(\mathbf{c})=\left(
			\begin{array}{cccc}
				\alpha_1 & 0 & 0 & \alpha_4 \\
				0& \beta_2 & 1-\alpha_1&0
			\end{array}\right)\simeq \left(
			\begin{array}{cccc}
				\alpha_1 & 0 & 0 & a^2\alpha_4 \\
				0& \beta_2 & 1-\alpha_1&0
			\end{array}\right),$\ where\\ $\mathbf{c}=(\alpha_1,\alpha_4, \beta_2)\in \mathbb{F}^3$ and $a\in \mathbb{F}^*$, \\
			\item	$ A_{4}(\mathbf{c})=\left(
			\begin{array}{cccc}
				0 & 1 & 1 & 0 \\
				\beta _1& \beta_2 & 1&-1
			\end{array}\right),\ \mbox{where}\ \mathbf{c}=(\beta_1, \beta_2)\in \mathbb{F}^2,$\\ 
			\item	$A_{5}(\mathbf{c})=\left(
			\begin{array}{cccc}
				\alpha _1 & 0 & 0 & 0 \\
				1 & 2\alpha _1-1 & 1-\alpha _1&0
			\end{array}\right),\ \mbox{where}\ \mathbf{c}=\alpha_1\in \mathbb{F},$\\
			\item	$ A_{6}(\mathbf{c})=\left(
			\begin{array}{cccc}
				\alpha_1 & 0 & 0 & \alpha_4 \\
				1& 1-\alpha_1 & -\alpha_1&0
			\end{array}\right),\ \mbox{where}\ \mathbf{c}=(\alpha_1,\alpha_4)\in \mathbb{F}^2,\ \alpha_4\in \mathbb{F}^*,$ \\ 
			\item	$ A_{7}(\mathbf{c})=\left(
			\begin{array}{cccc}
				\alpha_1 & 0 & 0 & \alpha_4 \\
				0&1-\alpha_1 & -\alpha_1&0
			\end{array}\right)\simeq \left(
			\begin{array}{cccc}
				\alpha_1 & 0 & 0 & a^2\alpha_4 \\
				0& 1-\alpha_1 & -\alpha_1&0
			\end{array}\right),$\ where\\ $\mathbf{c}=(\alpha_1,\alpha_4)\in \mathbb{F}^2$ and $ a\in \mathbb{F}^*$,\\
			\item	$ A_{8}(\mathbf{c})=\left(
			\begin{array}{cccc}
				0 & 1 & 1 & 0 \\
				\beta _1& 1 & 0&-1
			\end{array}\right),\ \mbox{where}\ \mathbf{c}=\beta_1\in\mathbb{F},$ \\
			\item	$ A_{9}(\mathbf{c})=\left(
			\begin{array}{cccc}
				\frac{1}{3} & 0 & 0 & 0 \\
				1 &\frac{2}{3} &-\frac{1}{3}&0
			\end{array}\right),$\\
			\item	$A_{10}(\mathbf{c})=\left(
			\begin{array}{cccc}
				0 &1 & 1 &1 \\
				\beta_1 &0 &0 &-1
			\end{array}
			\right)\simeq \left(
			\begin{array}{cccc}
				0 &1 & 1 &1 \\
				\beta'_1(a) &0 &0 &-1
			\end{array}\right),$ where $\mathbf{c}=\beta_1$,\\ $(\beta _1a^3-3a-1)(\beta _1a^2+\beta_1a+1)(\beta_1^2a^3+6\beta_1a^2+3\beta_1a+\beta_1-2)\neq 0,\\  \beta' _1(t)=\frac{(\beta_1^2t^3+6\beta_1t^2+3\beta_1t+\beta_1-2)^2}{(\beta_1t^2+\beta_1t+1)^3}$ if $t\neq \frac{-1}{2}$,  $\beta'(\frac{-1}{2})=4-\beta_1$.\\
			\item	$A_{11}(\mathbf{c})=\left(
			\begin{array}{cccc}
				0 &0 & 0 &1 \\
				\beta_1 &0 &0 &0
			\end{array}
			\right)\simeq \left(
			\begin{array}{cccc}
				0 &0 & 0 &1 \\
				a^3\beta_1^2 &0 &0 &0
			\end{array}
			\right)\simeq 
			\left(
			\begin{array}{cccc}
				0 &0 & 0 &1 \\
				b^3\beta_1 &0 &0 &0
			\end{array}\right),\ where\\ \ \mathbf{c}=\beta_1, a,b\in \mathbb{F}^*.$\\ Algebras $A_{10}(\beta_1)$ and $A_{11}(\beta'_1)$ are isomorphic if and only if there exists $t\in F^*$ such that $\beta_1=\beta'_1t^3+2+1/(\beta'_1t^3)$ and $(\beta'_1)^2t^6\neq 1$.\\
			\item	$A_{12}(\mathbf{c})=\left(
			\begin{array}{cccc}
				0 & 1 & 1 &0  \\
				\beta_1 &0& 0 &-1
			\end{array}
			\right)\simeq \left(
			\begin{array}{cccc}
				0 & 1 & 1 & 0 \\
				a^2\beta_1 &0& 0 &-1
			\end{array}
			\right),\ \mbox{where}\ \mathbf{c}=\beta_1\in \mathbb{F},\ a\in \mathbb{F}^*.$\\ Algebras $A_{10}(\beta_1)$ and $A_{12}(\beta'_1)$ are isomorphic  if and only if there exists $t\neq \pm1/2, s\neq 0$ such that $\beta_1=2(2t+1)^2(1-t), s^2\beta'_1=1-t^2$.\\ Algebras $A_{11}(\beta_1)$ and $A_{12}(\beta'_1)$ are isomorphic  if and only if there exists $s,t\in F^*$ such that $\beta_1=t^3, \beta'_1=-s^2$.\\
	\end{itemize}\end{theorem}
	\begin{theorem}Any non-trivial two-dimensional algebra over a field $\mathbb{F}$ with $\mathrm{Char}(\mathbb{F})=2$ is isomorphic to exactly one of the following algebras, listed by their matrices of structural constants:
		\begin{itemize}
			\item
			$ A_{1,2}(\mathbf{c})=\left(
			\begin{array}{cccc}
				\alpha_1 & \alpha_2 &\alpha_2+1 & \alpha_4 \\
				\beta_1 & \alpha_1 & 1+\alpha_1 & \alpha_2
			\end{array}\right),\ \mbox{where}\ \mathbf{c}=(\alpha_1, \alpha_2, \alpha_4, \beta_1) \in \mathbb{F}^4,$\\
			\item	$ A_{2,2}(\mathbf{c})=\left(
			\begin{array}{cccc}
				\alpha_1 & 0 & 0 & \alpha_4 \\
				1& \beta_2 & 1+\alpha_1&0
			\end{array}\right),\ \mbox{where}\ \mathbf{c}=(\alpha_1,\alpha_4, \beta_2)\in \mathbb{F}^3,\ \alpha_4\in \mathbb{F}^*,$\\
			$A_{2,2}(\alpha_1,0,1)=\left(
			\begin{array}{cccc}
				\alpha_1 & 0 & 0 &0  \\
				1& 1 & 1+\alpha_1&0
			\end{array}\right),\ \mbox{where}\ \alpha_1\in\mathbb{F}$,\\  
			\item	$ A_{3,2}(\mathbf{c})=\left(
			\begin{array}{cccc}
				\alpha_1 & 0 & 0 & \alpha_4 \\
				0& \beta_2 & 1+\alpha_1&0
			\end{array}\right)\simeq \left(
			\begin{array}{cccc}
				\alpha_1 & 0 & 0 & a^2\alpha_4 \\
				0& \beta_2 & 1+\alpha_1&0
			\end{array}\right),$ where\\ $\mathbf{c}=(\alpha_1,\alpha_4, \beta_2)\in \mathbb{F}^3$ and $a\in \mathbb{F}^*$, \\
			\item	$ A_{4,2}(\mathbf{c})=\left(
			\begin{array}{cccc}
				\alpha_1 & 1 & 1 & 0 \\
				\beta _1& \beta_2 & 1+\alpha_1&1
			\end{array}\right)\simeq \left(
			\begin{array}{cccc}
				\alpha_1 & 1 & 1 & 0 \\
				\beta _1+(1+\beta_2)a+a^2& \beta_2 & 1+\alpha_1&1
			\end{array}\right),\ \mbox{where}\\ \mathbf{c}=(\alpha_1,\beta_1, \beta_2)\in \mathbb{F}^2,$\\
			\item	$ A_{5,2}(\mathbf{c})=\left(
			\begin{array}{cccc}
				\alpha_1 & 0 & 0 & \alpha_4 \\
				1&1+\alpha_1 & \alpha_1&0
			\end{array}\right),\ \mbox{where}\ \mathbf{c}=(\alpha_1,\alpha_4)\in \mathbb{F}^2,\ \alpha_4\in \mathbb{F}^*,\\ A_{5,2}(1,0)=\left(
			\begin{array}{cccc}
				1 & 0 & 0 & 0 \\
				1&0 & 1&0
			\end{array}\right),$\\
			\item	$ A_{6,2}(\mathbf{c})=\left(
			\begin{array}{cccc}
				\alpha_1 & 0 & 0 & \alpha_4 \\
				0&1+\alpha_1 & \alpha_1&0
			\end{array}\right)\simeq \left(
			\begin{array}{cccc}
				\alpha_1 & 0 & 0 & a^2\alpha_4 \\
				0& 1+\alpha_1 & \alpha_1&0
			\end{array}\right),$ where\\ $\mathbf{c}=(\alpha_1,\alpha_4)\in \mathbb{F}^2$ and $ a\in \mathbb{F}^*$,\\
			\item	$ A_{7,2}(\mathbf{c})=\left(
			\begin{array}{cccc}
				\alpha_1 & 1 & 1 & 0 \\
				\beta _1& 1+\alpha_1 & \alpha_1&1
			\end{array}\right)\simeq \left(
			\begin{array}{cccc}
				\alpha_1 & 1 & 1 & 0 \\
				\beta _1+a\alpha_1+a+a^2& 1+\alpha_1 & \alpha_1&1
			\end{array}\right),\ \mbox{where}\\ \mathbf{c}=(\alpha_1,\beta_1)\in\mathbb{F}^2,\ a\in \mathbb{F},$\\
			\item	$ A_{8,2}(\mathbf{c})=\left(
			\begin{array}{cccc}
				0 &1 & 1 &1 \\
				\beta_1 &0 &0 &-1
			\end{array}
			\right)\simeq \left(
			\begin{array}{cccc}
				0 &1 & 1 &1 \\
				\beta'_1(a) &0 &0 &-1
			\end{array}\right),$\ where $\mathbf{c}=\beta_1$,\\ $(\beta _1a^3+a+1)(\beta _1a^2+\beta_1a+1)\beta_1\neq 0$,  and   $\beta' _1(t)=\beta_1^2\frac{(\beta_1t^3+t+1)^2}{(\beta_1t^2+\beta_1t+1)^3},$\\
			\item	$A_{9,2}(\mathbf{c})=\left(
			\begin{array}{cccc}
				0 &0 & 0 &1 \\
				\beta_1 &0 &0 &0
			\end{array}\right)
			\simeq \left(
			\begin{array}{cccc}
				0 &0 & 0 &1 \\
				a^3\beta_1^2 &0 &0 &0
			\end{array}
			\right)\simeq \left(
			\begin{array}{cccc}
				0 &0 & 0 &1 \\
				b^3\beta_1 &0 &0&0 
				\end{array}\right)$,
				 where  $\mathbf{c}=\beta_1\in \mathbb{F}$, $a,b\in \mathbb{F}^*$.\\
			Algebras $A_{8,2}(\beta_1)$ and $A_{9,2}(\beta'_1)$ are isomorphic if and only if there exists $t\in F^*$ such that $\beta_1=\beta'_1t^3+1/(\beta'_1t^3)$ and $(\beta'_1)^2t^6\neq 1$. 
			\item	$A_{10,2}(\mathbf{c})=\left(
			\begin{array}{cccc}
				1 & 1 & 1 & 0 \\
				\beta_1 &1& 1 &1
			\end{array}
			\right)\simeq \left(
			\begin{array}{cccc}
				1 & 1 & 1 & 0 \\
				\beta_1+a+a^2 &1& 1 &1
			\end{array}
			\right),\ \mbox{where}\ \mathbf{c}=\beta_1\in \mathbb{F},\ a\in \mathbb{F},$\\
			Algebras $A_{8,2}(\beta_1)$ and $A_{10,2}(\beta'_1)$ are isomorphic if and only if there exists $1\neq s\in F^*$ such that  $\beta_1=s^3+s$, $\beta'_1=(1+s)/s$.\\
			Algebras $A_{9,2}(\beta_1)$ and $A_{10,2}(\beta'_1)$ are isomorphic if and only if there exists $0\neq s,t\in F$ such that $\beta_1=s^3$, $t^2+st+(\beta'_1+1)s^2=0$. \item	$A_{11,2}(\mathbf{c})=\left(
			\begin{array}{cccc}
				0 & 1 & 1 & 0 \\
				\beta_1 &0& 0 &1
			\end{array}
			\right)\simeq \left(
			\begin{array}{cccc}
				0 & 1 & 1 & 0 \\
				b^2(\beta_1+a^2) &0& 0 &1
			\end{array}
			\right),$ where $\mathbf{c}=\beta_1$, $0\neq b$, $a,b\in \mathbb{F}$.\\
			Algebras $A_{8,2}(\beta_1)$ and $A_{11,2}(\beta'_1)$ are isomorphic if and only if  $\beta_1=0, \beta'_1=b^2$ for some $b\in \mathbb{F}$.\\
	\end{itemize} \end{theorem}
	\begin{theorem}Any non-trivial two-dimensional algebra over a field $\mathbb{F}$ with $\mathrm{Char}(\mathbb{F})=3$ is isomorphic to exactly one of the following algebras, listed by their matrices of structural constants:\begin{itemize}
			\item
			$ A_{1,3}(\mathbf{c})=\left(
			\begin{array}{cccc}
				\alpha_1 & \alpha_2 &\alpha_2+1 & \alpha_4 \\
				\beta_1 & -\alpha_1 & 1-\alpha_1 & -\alpha_2
			\end{array}\right),\ \mbox{where}\ \mathbf{c}=(\alpha_1, \alpha_2, \alpha_4, \beta_1) \in \mathbb{F}^4,$ 
			\item	$ A_{2,3}(\mathbf{c})=\left(
			\begin{array}{cccc}
				\alpha_1 & 0 & 0 & \alpha_4 \\
				1& \beta_2 & 1-\alpha_1&0
			\end{array}\right),\ \mbox{where}\ \mathbf{c}=(\alpha_1,\alpha_4, \beta_2)\in \mathbb{F}^3,\ \alpha_4\in \mathbb{F}^*,$
			\item	$ A_{3,3}(\mathbf{c})=\left(
			\begin{array}{cccc}
				\alpha_1 & 0 & 0 & \alpha_4 \\
				0& \beta_2 & 1-\alpha_1&0
			\end{array}\right)\simeq \left(
			\begin{array}{cccc}
				\alpha_1 & 0 & 0 & a^2\alpha_4 \\
				0& \beta_2 & 1-\alpha_1&0
			\end{array}\right),$ \ where\\ $\mathbf{c}=(\alpha_1,\alpha_4, \beta_2)\in \mathbb{F}^3$ and $a\in \mathbb{F}^*$, 
			\item	$ A_{4,3}(\mathbf{c})=\left(
			\begin{array}{cccc}
				0 & 1 & 1 & 0 \\
				\beta _1& \beta_2 & 1&-1
			\end{array}\right),\ \mbox{where}\ \mathbf{c}=(\beta_1, \beta_2)\in \mathbb{F}^2,$\
			\item	$ A_{5,3}(\mathbf{c})=\left(
			\begin{array}{cccc}
				\alpha _1 & 0 & 0 & 0 \\
				1 & 2\alpha _1-1 & 1-\alpha _1&0
			\end{array}\right),\ \mbox{where}\ \mathbf{c}=\alpha_1\in \mathbb{F},$
			\item	$ A_{6,3}(\mathbf{c})=\left(
			\begin{array}{cccc}
				\alpha_1 & 0 & 0 & \alpha_4 \\
				1& 1-\alpha_1 & -\alpha_1&0
			\end{array}\right),\ \mbox{where}\ \mathbf{c}=(\alpha_1,\alpha_4)\in \mathbb{F}^2,\ \alpha_4\in \mathbb{F}^*,$
			\item	$ A_{7,3}(\mathbf{c})=\left(
			\begin{array}{cccc}
				\alpha_1 & 0 & 0 & \alpha_4 \\
				0&1-\alpha_1 & -\alpha_1&0
			\end{array}\right)\simeq \left(
			\begin{array}{cccc}
				\alpha_1 & 0 & 0 & a^2\alpha_4 \\
				0& 1-\alpha_1 & -\alpha_1&0
			\end{array}\right),$ where\\ $\mathbf{c}=(\alpha_1,\alpha_4)\in \mathbb{F}^2$ and $a\in \mathbb{F}^*$,\ \
			\item	$ A_{8,3}(\mathbf{c})=\left(
			\begin{array}{cccc}
				0 & 1 & 1 & 0 \\
				\beta _1& 1 & 0&-1
			\end{array}\right),$\ where\ $\mathbf{c}=\beta_1\in\mathbb{F},$\\
			\item $A_{9,3}(\mathbf{c})=\left(
			\begin{array}{cccc}
				0 &1 & 1 &1 \\
				\beta_1 &0 &0 &-1
			\end{array}
			\right)\simeq \left(
			\begin{array}{cccc}
				0 &1 & 1 &1 \\
				\beta'_1(a) &0 &0 &-1
			\end{array}\right),$ where $\mathbf{c}=\beta_1$, $a\in \mathbb{F}$,\\ $(\beta _1a^3-1)(\beta _1a^2+\beta_1a+1)(\beta_1^2a^3+\beta_1+1)\neq 0$,  $\beta' _1(t)=\frac{(\beta_1^2t^3+\beta_1+1)^2}{(\beta_1t^2+\beta_1t+1)^3}$ if $t\neq 1$,  $\beta'(1)=2\beta_1+1$.
			
			\item $A_{10,3}(\mathbf{c})=\left(
			\begin{array}{cccc}
				0 &0 & 0 &1 \\
				\beta_1 &0 &0 &0
			\end{array}
			\right)\simeq \left(
			\begin{array}{cccc}
				0 &0 & 0 &1 \\
				a^3\beta_1^2 &0 &0 &0
			\end{array}
			\right)\simeq \left(
			\begin{array}{cccc}
				0 &0 & 0 &1 \\
				b^3\beta_1 &0 &0 &0
			\end{array}\right),$\ where  $\mathbf{c}=\beta_1$, $a,b\in \mathbb{F}^*$.\\
			Algebras $A_{9,3}(\beta_1)$ and $A_{10,3}(\beta'_1)$ are isomorphic if and only if there exists $t\in F^*$ such that $\beta_1=\beta'_1t^3+2+1/(\beta'_1t^3)$ and $(\beta'_1)^2t^6\neq 1$. \\
			\item $A_{11,3}(\mathbf{c})=\left(
			\begin{array}{cccc}
				0 &1 & 1 &0 \\
				\beta_1 &0 &0 &-1
			\end{array}
			\right)\simeq \left(
			\begin{array}{cccc}
				0 &1 & 1 &0 \\
				a^2\beta_1 &0 &0 &-1
			\end{array}\right),$\ where  $\mathbf{c}=\beta_1$, $a\in \mathbb{F}^*$.\\
			Algebras $A_{9,3}(\beta_1)$ and $A_{11,3}(\beta'_1)$ are isomorphic  if and only if there exists $\pm 1\neq t\in F, s\neq 0$ such that $\beta_1=t^3-1, s^2\beta'_1=1-t^2$. \\
			Algebras $A_{10,3}(\beta_1)$ and $A_{11,3}(\beta'_1)$ are isomorphic  if and only if there exists $s,t\in F^*$ such that $\beta_1=t^3, \beta'_1=-s^2$.
			
	\end{itemize}\end{theorem}
	
	\section{ Automorphisms and derivations of $2$-dimensional algebras}The following theorem describes the automorphisms and derivations, up to similarity, of all two-dimensional algebras in the case $\mathrm{Char}(\mathbb{F})\neq 2,3$.
	
	\begin{theorem}\label{thm4}
		The automorphisms and derivations of the algebras listed in Theorem 2.5 are given as follows.
		\begin{itemize}
			\item $Aut(A_1(\alpha_1,\alpha_2,\alpha_4,\beta_1))=\{I\},\\ Der(A_{1}(\alpha_1, \alpha_2, \alpha_4, \beta_1))=\{0\},$
			\item $Aut(A_2(\alpha_1, \alpha_4,\beta_2))=\{I\}$,\\
			$Der(A_{2}(\alpha_1,\alpha_4, \beta_2))= \{0\},$
			\item $Aut(A_3(\alpha_1,\alpha_4,\beta_2))=
			\{\left(
			\begin{array}{cc}
				1 & 0 \\
				0 & \pm1 
			\end{array}
			\right) \}$, if $\alpha_4\neq 0$,\\
						$Der(A_3(\alpha_1,\alpha_4,\beta_2))=\{ 0 \}$, if  $\alpha_4\neq 0$,\\
			$Aut(A_3(\alpha_1,0,\beta_2))=\{\left(
			\begin{array}{cc}
				1 & 0 \\
				0 & d 
			\end{array}
			\right); d\in \mathbb{F}, d\neq 0\},$ if $\beta_2\neq 2\alpha_1-1$,\\
			$Der(A_{3}(\alpha_1, 0, \beta_2))=\{\left(
			\begin{array}{cc}
				0 & 0 \\
				0 & d 
			\end{array}
			\right)
			:\ d\in \mathbb{F}\},$ if $ \beta_2 \neq 2 \alpha_1-1,$\\
			$Aut(A_3(\alpha_1,0,2\alpha_1-1))=\{\left(
			\begin{array}{cc}
				1 & 0 \\
				c & d 
			\end{array}
			\right); c,d\in \mathbb{F}, d\neq 0\},$\\
			$Der(A_{3}(\alpha_1, 0, 2\alpha_1-1))=\{\left(
			\begin{array}{cc}
				0 & 0 \\
				c & d
				\end{array}\right):\ c, d\in \mathbb{F}\}$, 
			\item $Aut(A_4(\beta_1,\beta_2))=\{I\},$\\
			$Der(A_{4}(\beta_1, \beta_2))=\{ 0\},$ 
			\item $Aut(A_5(\alpha_1))=\{ \left(
			\begin{array}{cc}
				1 & 0 \\
				c & 1\end{array}\right):c \in\mathbb{F}\},$\\
			$Der(A_{5}(\alpha_1))=\{\left(
			\begin{array}{cc}
				0 & 0 \\
				c & 0\end{array} \right):\ c\in \mathbb{F}\},$
			\item $Aut(A_6(\alpha_1,\alpha_4))=\{I\},$\\
			$Der(A_{6}(\alpha_1, \beta_1))=\{0\},$
			\item 
			$Aut(A_7(\alpha_1, \alpha_4))=\{\begin{pmatrix}
				1 & 0 \\ 
				0& \pm 1
			\end{pmatrix}\},$ if $\alpha_4\neq 0$,\\
			$Der(A_7(\alpha_1, \alpha_4))=\{0\},$ if $\alpha_4 \neq 0$,\\ $Aut(A_7(\alpha_1,0))=\{\begin{pmatrix}
				1 & 0 \\ 
				0& d
			\end{pmatrix}: 0\neq d\in\mathbb{F} \}$ if $\alpha_1\neq \frac{1}{3}$,\\
			$Der(A_{7}(\alpha_1, 0))=\{\left(
			\begin{array}{cc}
				0 & 0 \\
				0 & d\end{array}\right)
			:\ d\in \mathbb{F}\},$ if $\alpha_1 \neq \frac{1}{3},$\\ $Aut(A_7(\frac{1}{3},0))=\{\begin{pmatrix}
				1 & 0 \\ 
				c& d
			\end{pmatrix}: 0\neq d,c\in\mathbb{F}\},$\\
			$	Der(A_{7}(\frac{1}{3}, 0))=\{\left(
			\begin{array}{cc}
				0 & 0 \\
				c & d 
				\end{array}\right): c, d\in \mathbb{F}\},$
			\item $Aut(A_8(\beta_1))=\{I\},$\\
			$Der(A_{8}(\beta_1))=\{0\},$
			\item $Aut(A_9)=\{\left(
			\begin{array}{cc}
				1 & 0 \\
				c & 1\end{array}\right):c \in\mathbb{F}\},$\\
			$ Der(A_{9})=\{\left(
			\begin{array}{cc}
				0 & 0 \\
				c & 0\end{array} \right):\ c \in \mathbb{F}\},$
			\item$Aut(A_{10}(\beta_1))=\\ \left\{I\right\}\cup \{\begin{pmatrix}
				-1-d&-\frac{1+d+d^2}{1+2d}\\ 1+2d&d
			\end{pmatrix}:\ \beta_1=\frac{(1+2d)^2}{d^2+d+1}\}\cup \{\begin{pmatrix}
				-d & \frac{-d-d^2}{1+2d}\\
				\frac{2d^2-d-1}{d} & d 
			\end{pmatrix}:\ \beta_1=\frac{(d-1)(1+2d)^2}{d^3}\}$, for $\beta_1\neq 0,4$, \\
			$Der(A_{10}(\beta_1))=\{0\},$ if $\beta_1\neq 0,4$,\\
			$Aut(A_{10}(0))=\{\left(
			\begin{array}{cc}
				3a+1& a \\
				0 & 1 
			\end{array}\right): \ \ 3a+1\neq 0\}$,\\
			$Der(A_{10}(0)=\{\begin{pmatrix}
				3a & a\\
				0 & 0
			\end{pmatrix}:\ a\in\mathbb{F}\}$,\\
			\item $Aut(A_{11}(\beta_1))=\{\left(
			\begin{array}{cc}
				a & 0 \\
				0 & a^2 
			\end{array}
			\right):\ a^3=1\}\cup \{\begin{pmatrix}0&b\\  
				1/b&0\end{pmatrix}:\ b^{-3}=\beta_1\},$ if $\beta_1\neq 0$ ,\\
					$Der(A_{11}(\beta_1))=\{0\}$ if $\beta_1\neq 0$,\\
			$Aut(A_{11}(0))=\{\left(
			\begin{array}{cc}
				a^2 & b \\
				0 & a 
			\end{array}
			\right):\ a,b\in \mathbb{F},\ a\neq 0\},$\\
					$Der(A_{11}(0))=\{\begin{pmatrix}
				2a & b\\
				0 & a
			\end{pmatrix}: \ a,b\in \mathbb{F}\},$   
						\item $Aut(A_{12}(\beta_1))=\{\left(
			\begin{array}{cc}
				\pm1 & 0 \\
				0 & 1 
			\end{array}
			\right)\}\cup \{\left(
			\begin{array}{cc}
				\pm\frac{1}{2} & b \\
				\pm\frac{3}{4b} & -\frac{1}{2}  
			\end{array}
			\right):\ \beta_1=\frac{3}{4b^2}, b\in\mathbb{F} \},$ if $\beta_1\neq 0$,\\ 
			$Der(A_{12}(\beta_1))=\{0\},$ if $\beta_1 \neq 0,$\\
			$Aut(A_{12}(0))=\{\left(
			\begin{array}{cc}
				a & 0 \\
				0 & 1 
			\end{array}
			\right):\ 0\neq a\in\mathbb{F} \},$ \\
			$Der(A_{12}(0))=\{\begin{pmatrix}
				a & 0\\
				0 & 0
			\end{pmatrix}:a \in \mathbb{F}\},$\\
			
		\end{itemize}
	\end{theorem}
	\textbf{Proof.} Let us first prove the automorphism part of the theorem.
	
	To prove the theorem, we solve equation \eqref{1} for $g\in GL(2,\mathbb{F})$ in each case $A=A_i$, where $i=1,2,\ldots,12$.
	Note that, by \cite{ABR}, if $A=\left(
	\begin{array}{cccc}
		\alpha_1 & \alpha_2 &\alpha_3 & \alpha_4 \\
		\beta_1 & \beta_2 & \beta_3 & \beta_4
	\end{array}\right)$ and $P(A)=\left(
	\begin{array}{c}
		\overline{tr}_1(A)\\ \overline{tr}_2(A)
	\end{array}\right)$, where $$\overline{tr}_1(A)=(\alpha_1+\beta_3, \alpha_2+\beta_4), \overline{tr}_2(A)=(\alpha_1+\beta_2, \alpha_3+\beta_4)$$-are row vectors, then for every $g\in GL(2,\mathbb{F})$ the equality $P(gA(g^{-1})^{\otimes 2})=P(A)g^{-1}$ holds. 
	
	Therefore, in the case
	$ A=A_{1}(\alpha_1, \alpha_2, \alpha_4, \beta_1)=\left(
	\begin{array}{cccc}
		\alpha_1 & \alpha_2 &\alpha_2+1 & \alpha_4 \\
		\beta_1 & -\alpha_1 & -\alpha_1+1 & -\alpha_2
	\end{array}\right)$ the equality $A=g^{-1}Ag^{\otimes 2}$ holds  if and only if $g=I$, since in this case $P(A)=I$. Hence	
	$Aut(A_1(\alpha_1,\alpha_2,\alpha_4,\beta_1))=\{I\}.$ This proof does not depend on $Char(\mathbb{F})$.
	
	Note that in the cases $A=A_{i}(c)$, $i=2,3,4,5$, one has   $\overline{tr}_1(A)=(1,0)$ and in the cases $A=A_{i}(c)$, $i=6,7,8,9$, one has $\overline{tr}_2(A)=(1,0)$. Therefore, in these cases, if
	$A=g^{-1}Ag^{\otimes 2}$, then $(1,0)=(1,0)g$, which implies that $g=\left(
	\begin{array}{cc}
		1& 0 \\
		c & d 
	\end{array}
	\right)$, where $d\neq 0$. Therefore, in all these cases we consider the equality $gA=Ag^{\otimes 2}$ for matrices of this form. 
	
	If $A=A_{2}(\mathbf{c})=\left(
	\begin{array}{cccc}
		\alpha_1 & 0 & 0 & \alpha_4 \\
		1& \beta_2 & 1-\alpha_1&0
	\end{array}\right),\ \mbox{where}\ \mathbf{c}=(\alpha_1,\alpha_4, \beta_2)\in \mathbb{F}^3,\ \alpha_4\neq 0$, then the equality $gA=Ag^{\otimes 2}$ is equivalent to
	\[\begin{array}{lc}
	c^2 \alpha_4  &=0,\\
	-c d \alpha_4 &=0,\\
	\alpha_4 - d^2 \alpha_4 &=0,\\
	- c - 1+d+2c \alpha _1- c \beta _2 &=0,\\
	c \alpha_4&=0.\end{array}\] 
	Hence $c=0$ and $d=1$ because of $\alpha_4 \neq 0$. Therefore, $Aut(A_2(c))=\{I\}$.
	
	If $ A_{3}(\mathbf{c})=\left(
	\begin{array}{cccc}
		\alpha_1 & 0 & 0 & \alpha_4 \\
		0& \beta_2 & 1-\alpha_1&0
	\end{array}\right)\simeq \left(
	\begin{array}{cccc}
		\alpha_1 & 0 & 0 & a^2\alpha_4 \\
		0& \beta_2 & 1-\alpha_1&0
	\end{array}\right),$ where $\mathbf{c}=(\alpha_1,\alpha_4, \beta_2)\in \mathbb{F}^3$ and $0\neq a\in \mathbb{F}$, then $gA=Ag^{\otimes 2}$ is equivalent to
	\[\begin{array}{lc}
	c \alpha_4& =0,\\
	\alpha_4 - d^2 \alpha_4&  =0,\\
	-c+2c \alpha _1- c \beta _2 &=0.\end{array}\]
	If $\alpha_4=0$, then the system becomes 
	$
	c(-1+2\alpha _1-\beta _2) =0$, 
	hence $Aut(A_{3}(\alpha_1,0, \beta_2))=\{\left(\begin{array}{cc}
		1& 0\\
		0&d
	\end{array}\right): d\neq 0\}$, if $\beta _2\neq 2\alpha _1-1$,
	$Aut(A_{3}(\alpha_1,0, 2\alpha _1-1))=\{\left(\begin{array}{cc}
		1& 0\\
		c&d
	\end{array}\right): d\neq 0\}$.	
	If $\alpha_4\neq 0$, then $c=0, d=\pm 1$ and $Aut(A_{3}(\alpha_1,\alpha_4, \beta_2))=\{\left(\begin{array}{cc}
		1& 0\\
		0&\pm 1
	\end{array}\right)\}$. 
	
	In the case $ A=A_{4}(\beta_1, \beta_2)=\left(
	\begin{array}{cccc}
		0 & 1 & 1 & 0 \\
		\beta_1& \beta _2 & 1&-1
	\end{array}\right), $ we obtain the system
\[\begin{array}{lc}		-2 c&=0,\\
	1-d&=0,\\
	-c+c^2-\beta _1+d \beta _1- c \beta _2& =0,\\
	c+c d& =0,\\
	c+c d&=0,\\
	-d+d^2&=0,\end{array}\] hence $c=0,d=1$ and
	$Aut(A_{4}c))=\{I\}$.
	
	If $ A=A_{5}(\alpha_1)=\left(
	\begin{array}{cccc}
		\alpha_1& 0 & 0 & 0 \\
		1 & 2\alpha_1-1 & 1-\alpha_1&0
	\end{array}\right)$, then the system becomes
	$		 d-1=0$, hence  $Aut(A_5(\alpha_1))=\{ \left(
	\begin{array}{cc}
		1 & 0 \\
		c & 1 
	\end{array}
	\right):c \in\mathbb{F}\}.$
	
	In the case $ A_{6}(\alpha_1, \alpha_4)=\left(
	\begin{array}{cccc}
		\alpha_1 & 0 & 0 & \alpha_4 \\
		1& 1-\alpha_1 & -\alpha_1&0
	\end{array}\right), \alpha_4 \neq 0,$ the system becomes \[\begin{array}{lc}
	c \alpha_4 &=0,\\
	\alpha_4 - d^2 \alpha_4 &=0,\\
	d-1 -c +3c \alpha_1&=0, 
\end{array}\]
	hence $c=0, d=1$ and $Aut(A_6(c))=\left\{I\right\}.$
	
	In the case $ A_{7}(\alpha_1, \alpha_4)=\left(
	\begin{array}{cccc}
		\alpha_1 & 0 & 0 & \alpha_4 \\
		0& 1-\alpha_1 & -\alpha_1&0
	\end{array}\right) \simeq \left(
	\begin{array}{cccc}
		\alpha_1 & 0 & 0 & t^2 \alpha_4 \\
		0& 1-\alpha_1 & -\alpha_1&0
	\end{array}\right)$ we obtain the system \[\begin{array}{lc}
	c\alpha_4&=0,\\
	\alpha_4 - d^2 \alpha_4& =0,\\
	3c\alpha_1-c&=0,\end{array}\]
	hence, if $\alpha_4\neq 0$, then $Aut(A_7(\alpha_1, \alpha_4))=\{\begin{pmatrix}
		1 & 0 \\ 
		0& \pm 1
	\end{pmatrix}\},$ otherwise $Aut(A_7(\alpha_1,0))=\{\begin{pmatrix}
		1 & 0 \\ 
		0& d
	\end{pmatrix}: 0\neq d\in\mathbf{F}\}$ if $\alpha_1\neq \frac{1}{3}$, and $Aut(A_7(\frac{1}{3},0))=\{\begin{pmatrix}
		1 & 0 \\ 
		c& d
	\end{pmatrix}: 0\neq d,c\in\mathbb{F}\}$. 
	
	In  the case $ A_{8}(\beta_1)=\left(
	\begin{array}{cccc}
		0 & 1 & 1 & 0 \\
		\beta_1& 1& 0&-1
	\end{array}\right)$ we obtain  
	$gA_{8}-A_{8}(g\otimes g)=$ \[\left(
	\begin{array}{cccc}
		-2c& 1-d &
		1-d &0\\
		-c+c^2-\beta _1+d \beta _1 & c+c d& c+c d& -d+d^2
	\end{array}
	\right),\] hence $Aut( A_{8}(\beta_1))=\left\{I\right\}$.
	
	In the case $ A=A_{9}=\left(
	\begin{array}{cccc}
		\frac{1}{3}& 0 & 0 & 0 \\
		1 & \frac{2}{3} & -\frac{1}{3}&0
	\end{array}\right)$ we have 
	$gA_{9}-A_{9}(g\otimes g)=\left(
	\begin{array}{cccc}
		0 &0&0&0\\
		-1+d & 0&0&0
	\end{array}
	\right),$ hence $Aut( A_{9})=\{\left(
	\begin{array}{cc}
		1 & 0 \\
		c & 1 
	\end{array}
	\right): c\in\mathbb{F}\}$.
	
	In the remaining cases we have to consider $gA=Ag^{\otimes 2}$ with respect to $g=\left(
	\begin{array}{cc}
		a& b \\
		c & d 
	\end{array}
	\right)$,\\ where $\Delta=ad-bc\neq 0$.
	
	Therefore, in the case $A=A_{10}(\beta_1)=\left(
	\begin{array}{cccc}
		0 & 1 & 1 & 1 \\
		\beta_1& 0& 0&-1
	\end{array}\right)\simeq \left(
	\begin{array}{cccc}
		0 & 1 & 1 & 1 \\
		4-\beta_1& 0& 0&-1
	\end{array}\right)=A_{10}(4-\beta_1),\ \mbox{where}\ \mathbf{c}=\beta_1\in \mathbb{F}$, the system becomes \[\begin{array}{lc}
	(1)\ \ -2ac -c^2+b\beta_1& = 0,\\
	(2)\ \ a-bc-ad-cd& = 0,\\
	(3)\ \	a-b-2bd-d^2&= 0,\\
	(4)\ \	c^2-a^2\beta_1+d\beta_1&= 0,\\
	(5)\ \	c+cd-ab\beta_1&= 0,\\
	(6)\ \	c-d+d^2-b^2\beta_1&= 0.\end{array}\]
	
	Equations (3) and (6) imply $a = b+2bd+d^2$ and $c = d-d^2+b^2\beta_1$, $\Delta=d^3+3bd^2-b^3\beta_1$. Substituting these expressions into the remaining equations, we obtain\\
	$d(1-d)(2b+d+4bd+d^2)+b(-1+2bd+b^2(2+4d))\beta_1+b^4\beta_1^2= 0,\\
	b(-1+d^2+b(b+d)\beta_1) = 0,\\
	d\beta_1-(b+2bd+d^2)^2\beta_1+(d-d^2+b^2\beta_1)^2 = 0,\\
	d(-1+d^2+b(b+d)\beta_1) = 0.$ 
	
	If $b=0$, then \[\begin{array}{lc} d(1-d)(d+d^2)&= 0,\\
	d\beta_1-d^4\beta_1+(d-d^2)^2& = 0,\\
	d(-1+d^2)& = 0.\end{array}\]
	Therefore, if $d=1$, then $g=I$; if $d=-1$, then $\beta_1=2$ and $g(2)=\left(
	\begin{array}{cc}
		1 & 0 \\
		-2 & -1 
	\end{array}
	\right).$
	
	In the case $b=-d$ the system becomes \[\begin{array}{lc} d(1-d)(-d-3d^2)-d(-1-2d^2+d^2(2+4d))\beta_1+d^4\beta_1^2&= 0,\\
	d(d^2-1) &= 0,\\
	d\beta_1-(d+d^2)^2\beta_1+(d-d^2+d^2\beta_1)^2 &= 0,
	\\	d(d^2-1) &= 0.\end{array}\]	Therefore, if $d=1$, then $\beta_1=0$ with $g(0)=\left(
	\begin{array}{cc}
		-2 & -1 \\
		0 & 1 
	\end{array}
	\right)$, or $\beta_1=3$ with $g(3)=\left(
	\begin{array}{cc}
		-2 & -1 \\
		3 & 1 
	\end{array}
	\right).$
	
	 If $d=-1$, then $\beta_1=1$ with $g(1)=\left(
	\begin{array}{cc}
		0 & 1 \\
		-1 & -1 
	\end{array}
	\right)$ or $\beta_1=4$, which can be ignored due to $A_{10}(0)\simeq A_{10}(4)$.
	
	If $b(b+d) \neq0$, then $\beta_1= \frac{1-d^2}{b(b+d)}$ and the system becomes \\ $\left\{\begin{array}{c}
		\frac{(b+d)(-1+d^2)+(-1+d^2)(b+2bd+d^2)^2-2(-1+d)(b+d)(b+2bd+d^2)^2}{(b+d)^2} =0\\
		\frac{(-1+d)^2(b+2bd+d^2)^2}{(b+d)^2}+\frac{(-1+d^2)(b+2bd+d^2)^2}{b(b+d)}+\frac{d(1-d^2)}{b(b+d)} =0\end{array}\right. \simeq \\ \left\{\begin{array}{c}
		(b+d)(1+d)-(1+2b+d)(b+2bd+d^2)^2=0\\
		d(1+2b+d)(b+2bd+d^2)^2=d(d+1)(b+d)\end{array}\right. \simeq$ \\ $(b+d)(1+d)-(1+2b+d)(b+2bd+d^2)^2=0$ and 
	$(b+d)(1+d)-(1+2b+d)(b+2bd+d^2)^2=$\\ $-2(1+2d)^2(b-\frac{1-d}{2})(b+\frac{1+d+d^2}{1+2d})(b+\frac{d+d^2}{1+2d}),$ provided that $1+2d\neq 0$.
	
	\underline{\textbf{Case 1.}} If  $b_1(d)=b= \frac{1-d}{2}$, then $\beta_1=4$ with $g(4)=\left(
	\begin{array}{cc}
		(1+d)/2 & (1-d)/2 \\
		1-d & d 
	\end{array}
	\right)$. Using $gA_{10}(0)(g^{-1})^{\otimes 2}=A_{10}(4)$, where  $g=\left(
	\begin{array}{cc}
		1 & 0 \\
		-2 & -1 
	\end{array}
	\right)$ we can conclude that $Aut(A_{10}(0))=\{\left(
	\begin{array}{cc}
		3a+1& a \\
		0 & 1 
	\end{array}\right): \ \ a\in \mathbb{F}, 3a+1\neq 0\}$.
	
		\underline{\textbf{Case 2.}}
	If $b_2(d)=b= \frac{-1-d-d^2}{1+2d}$, then $a=-1-d$, $c=1+2d$ and $g=\begin{pmatrix}
		-1-d & \frac{-1-d-d^2}{1+2d} \\
		1+2d & d
	\end{pmatrix}$, where $\beta_1=\beta_1(d)= \frac{(1+2d)^2}{d^2+d+1}$, $d^2+d+1 \neq 0$. 	
	Note that $\beta_1(-1-d)=\frac{(1+2d)^2}{d^2+d+1}=\beta_1(d)$.
	
	\underline{\textbf{Case 3.}} If $b_3(d)=b=\frac{-d-d^2}{1+2d}$, then $a=-d$, $c=\frac{d^2-d-1}{d}$ and $\Delta=-1.$ Therefore, $g=\begin{pmatrix}
		-d & \frac{-d-d^2}{1+2d}\\
		-1-\frac{1}{d}+2d & d 
	\end{pmatrix},$ where $\beta_1=\frac{(d-1)(1+2d)^2}{d^3}$. 
	
	If $1+2d=0$, then the last equation implies $b=3/4$, so $\beta_1= 4$, which can be omitted. Hence	
	$Aut(A_{10}(\beta_1))=\\ \left\{I\right\}\cup \{\begin{pmatrix}
		-1-d&-\frac{1+d+d^2}{1+2d}\\ 1+2d&d
	\end{pmatrix}:\ \beta_1=\frac{(1+2d)^2}{d^2+d+1}\}\cup \{\begin{pmatrix}
		-d & \frac{-d-d^2}{1+2d}\\
		\frac{2d^2-d-1}{d} & d 
	\end{pmatrix}:\ \beta_1=\frac{(d-1)(1+2d)^2}{d^3}\}$ for all $\beta_1\neq 0,4$. $Aut(A_{10}(0))=\{\left(
	\begin{array}{cc}
		3a+1& a \\
		0 & 1 
	\end{array}\right): \ \ a\in \mathbb{F}, 3a+1\neq 0\}$.

	If $ A=A_{11}(\beta_1)=\left(
	\begin{array}{cccc}
		0 & 0 & 0 & 1 \\
		\beta_1& 0 & 0& 0
	\end{array}\right)$,  then 
	$0=gA-A(g\otimes g)=$\\ $\left(
	\begin{array}{cccc}
		-c^2 +b \beta _1 & -c d &
		-c d & a-d^2 \\
		-a^2 \beta _1+d \beta _1 & -a b \beta _1 & -a b \beta _1 & c-b^2 \beta _1
	\end{array}
	\right).$ Hence, if $\beta _1=0$, then $c=0, a=d^2$ and $\begin{pmatrix}d^2&b\\  
		0&d\end{pmatrix}\in Aut(A_{11}(0))$, where $b,d\in \mathbb{F}, d\neq 0$. If $\beta_1\neq 0$, then $b=c=0$ or $a=d=0$ and $\begin{pmatrix}a&0\\  
		0&a^2\end{pmatrix}\in Aut(A_{11}(\beta_1))$, where $\beta_1\neq 0$, $a^3=1$, $\begin{pmatrix}0&b\\  
		1/b&0\end{pmatrix}\in Aut(A_{11}(\beta_1))$, where $b^{-3}=\beta_1$, respectively. It implies that\\
	$Aut(A_{11}(\beta_1))=\{\left(
	\begin{array}{cc}
		a & 0 \\
		0 & a^2 
	\end{array}
	\right):\ a^3=1\}\cup \{\begin{pmatrix}0&b\\  
		1/b&0\end{pmatrix}:\ b^{-3}=\beta_1\},$ in $\beta_1\neq 0$ case,\\
	$Aut(A_{11}(0))=\{\left(
	\begin{array}{cc}
		d^2 & b \\
		0 & d 
	\end{array}
	\right):\ b,d\in \mathbb{F},\ d\neq 0\}.$

	If $ A_{12}(\beta_1)=\left(
	\begin{array}{cccc}
		0 & 1 & 1 & 0 \\
		\beta_1& 0& 0&-1
	\end{array}\right),$  then
	$gA-A(g\otimes g)=0$ is equivalent to \[\begin{array}{lc}
	-2 a c+b \beta _1 &=0,\\
	a-b c-a d&=0,\\
	-b(1+2d) & =0,\\
	c^2-a^2 \beta _1+d \beta _1 &=0,\\
	c+c d-a b \beta _1 &=0, \\
	-d+d^2-b^2 \beta _1 &=0.\end{array}\]
	If $b=0$, then $c=0,d=1$ and $(1-a^2) \beta_1=0$. Hence  $g=\left(
	\begin{array}{cc}
		a & 0 \\
		0 & 1 
	\end{array}
	\right)$ if $\beta_1=0$ and $g=\left(
	\begin{array}{cc}
		\pm 1 & 0 \\
		0 & 1 
	\end{array}
	\right)$ \
	if  $\beta_1 \neq 0$.\\	
	If $b\neq 0$, then $d=\frac{-1}{2}$ and the system becomes \[\begin{array}{lc}
	-2 a c+b \beta _1 &=0,\\
	\frac{3a}{2} - b c&=0,\\
	c^2- \frac{1}{2} (1 + 2 a^2) \beta _1 &=0,\\
	\frac{c - 2 a b \beta_1}{2}&=0, \\
	\frac{3}{4} - b^2 \beta_1 &=0. \end{array}\]
	Hence $a=\frac{2 b c}{3}$ and
	\[\begin{array}{lc}
		b( \frac{-4 c^2}{3} + \beta _1) &=0,\\
		c^2 - \frac{9+ 8 b^2 c^2}{18} \beta_1 &=0,\\
		\frac{c}{6} (3- 4 b^2 \beta_1) &=0,\\
		\frac{3}{4} - b^2 \beta_1 &=0.\end{array}\]
	Therefore, $\beta_1=\frac{4 c^2}{3}$ and	
	\[\begin{array}{lc}
		\frac{c^2}{27}(9 - 16 b^2 c^2)&=0,\\
		\frac{c}{2} - \frac{8 b^2 c^3}{9} &=0,\\
		\frac{3}{4} - \frac{4 b^2 c^2}{3} &=0,\end{array}\]
	$c= \pm \frac{3}{4b}$. It means that $Aut(A_{12}(0))=\{\left(
	\begin{array}{cc}
		a & 0 \\
		0 & 1 
	\end{array}
	\right):\ 0\neq a\in\mathbb{F}\},$ \\ 
	$Aut(A_{12}(\beta_1))=\{\left(
	\begin{array}{cc}
		\pm1 & 0 \\
		0 & 1 
	\end{array}
	\right)\}\cup \{\left(
	\begin{array}{cc}
		\pm\frac{1}{2} & b \\
		\pm\frac{3}{4b} & -\frac{1}{2}  
	\end{array}
	\right):\ \beta_1=\frac{3}{4b^2}, b\in\mathbb{F} \},$ if $\beta_1\neq 0$.
	
		Here we prove the derivation part of the theorem.
		
	If $\mathbf{D}: \mathbb{A}\rightarrow\mathbb{A}$ is a derivation  with matrix $ D=\left(\begin{array}{cc}
		a &b \\
		c &d
	\end{array}\right)$,  then the equality
	$DA=A(D\otimes I+I\otimes D)$ is equivalent to 
	\begin{equation} \label{SE4}
		\begin{array}{lc}
			a \alpha _1-b \beta _1+\alpha _2 c+\alpha _3 c &=0,\\
			\alpha _1 b-b \beta _2+\alpha _4 c+\alpha _2 d &=0, \\
			\alpha _1 b-b \beta _3+\alpha _4 c+\alpha _3 d &=0,\\
			-a \alpha _4+\alpha _2 b+\alpha _3 b-b \beta _4+2 \alpha _4 d &=0,\\
			\beta _1 (2 a-d)-\alpha _1 c+\left(\beta _2+\beta _3\right) c &=0,\\
			a \beta _2+b \beta _1-\alpha _2 c+\beta _4 c &=0,\\
			a \beta _3+b \beta _1-\alpha _3 c+\beta _4 c &=0,\\
			b \beta _2+b \beta _3-\alpha _4 c+\beta _4 d &=0.\end{array}\end{equation}
	Moreover, the identities $\overline{tr}_1(A)D=\overline{tr}_2(A)D=0$ hold, where   $\overline{tr}_1(A)=(\alpha_1+\beta_3, \alpha_2+\beta_4), \overline{tr}_2(A)=(\alpha_1+\beta_2, \alpha_3+\beta_4)$.
	
	If $ A=A_{1}(\alpha_1, \alpha_2, \alpha_4, \beta_1)=\left(
	\begin{array}{cccc}
		\alpha_1 & \alpha_2 &\alpha_2+1 & \alpha_4 \\
		\beta_1 & -\alpha_1 & -\alpha_1+1 & -\alpha_2
	\end{array}\right)$, then $P(A)=I_2=I$, and therefore the equality $P(A)D=0$ implies $D=0.$
	
	If $\overline{tr}_1(A)=(1,0)$ or $\overline{tr}_2(A)=(1,0)$, which happens in $A_{2}(\mathbf{c}), A_{3}(\mathbf{c}), A_{4}(\mathbf{c}), A_{5}(\mathbf{c}), A_{6}(\mathbf{c}), A_{7}(\mathbf{c}), A_{8}(\mathbf{c})$ cases, the identity $P(A)D=0$ implies that $a=b=0$, that is $ D=\left(\begin{array}{cc}
		0 &0 \\
		c &d
	\end{array}\right)$.
	
	In the case $A=A_{2}(\mathbf{c})=\left(
	\begin{array}{cccc}
		\alpha_1 & 0 & 0 & \alpha_4 \\
		1& \beta_2 & 1-\alpha_1&0
	\end{array}\right)$ the system (\ref{SE4})  becomes
	\[\begin{array}{lc}
		c \alpha_4 &=0,\\
		2 d \alpha_4 &=0,\\
		c-d-2 c \alpha_1+c \beta_2  &=0.\end{array}\]
	Hence, since $\alpha_4 \neq 0$, we obtain $D=0$
	
	In the case $A=A_{3}(\mathbf{c})=\left(
	\begin{array}{cccc}
		\alpha_1 & 0 & 0 & \alpha_4 \\
		0& \beta_2 & 1-\alpha_1&0
	\end{array}\right)$ the system is 
	\[\begin{array}{lc}
		c \alpha_4 &=0,\\
		2 d \alpha_4 &=0,\\
		c(1-2 \alpha_1+\beta_2)  &=0.\end{array}\]
	Hence, if $\alpha_4=0$ and $\beta_2=2 \alpha_1-1$, then $D=\begin{pmatrix}
		0 & 0\\
		c & d
	\end{pmatrix};$\
	if $\alpha_4=0$ and $\beta_2 \neq 2 \alpha_1-1$, then $D=\begin{pmatrix}
		0 & 0\\
		0 & d
	\end{pmatrix};$\
	if $\alpha_4 \neq 0$, then $D=0.$
	
	If $A=A_{4}(\mathbf{c})=\left(
	\begin{array}{cccc}
		0 & 1 & 1 & 0 \\
		\beta _1& \beta_2 & 1&-1
	\end{array}\right)$, then the system becomes
	\[\begin{array}{lc}
		2 c  &=0,\\
		d &=0,\\
		c-d \beta_1+c \beta_2 &=0.\end{array}\] This implies that  $D=0.$
	
	In the case $A=A_{5}(\mathbf{c})=\left(
	\begin{array}{cccc}
		\alpha _1 & 0 & 0 & 0 \\
		1 & 2\alpha _1-1 & 1-\alpha _1&0
	\end{array}\right)$ the system (\ref{SE4}) is equivalent to the equation $d=0.$ Hence 
	$D=\begin{pmatrix}
		0 & 0\\
		c & 0
	\end{pmatrix}.$
	
	If $A=A_{6}(\mathbf{c})=\left(
	\begin{array}{cccc}
		\alpha_1 & 0 & 0 & \alpha_4 \\
		1& 1-\alpha_1 & -\alpha_1&0
	\end{array}\right)$ the system is
	\[\begin{array}{lc}
		c \alpha_4 &=0,\\
		2 d \alpha_4 &=0,\\
		c-d-3 c \alpha_1  &=0.\end{array}\] hence, due to $\alpha_4 \neq 0$, we have $D=0.$
	
	If $A=A_{7}(\mathbf{c})=\left(
	\begin{array}{cccc}
		\alpha_1 & 0 & 0 & \alpha_4 \\
		0&1-\alpha_1 & -\alpha_1&0
	\end{array}\right)$ the system is
	\[\begin{array}{lc}
		c \alpha_4 &=0,\\
		2 d \alpha _4  &=0,\\
		c-3 c \alpha _1  &=0,\end{array}\] Hence, if $\alpha_4=0$ and $\alpha_1=\frac{1}{3}$, then $D=\begin{pmatrix}
		0 & 0\\
		c & d
	\end{pmatrix};$ if $\alpha_4=0$ and $\alpha_1\neq\frac{1}{3}$, then
	$D=\begin{pmatrix}
		0 & 0\\
		0 & d
	\end{pmatrix};$ otherwise $D=0.$
	
	In the case $A=A_{8}(\mathbf{c})=\left(
	\begin{array}{cccc}
		0 & 1 & 1 & 0 \\
		\beta _1& 1 & 0&-1
	\end{array}\right)$ the system becomes
	\[\begin{array}{lc}
		2 c &=0,\\
		d  &=0,\\
		c-d \beta_1 &=0, \end{array}\]
	hence $D=0.$
	
	If $A=A_{9}(\mathbf{c})=\left(
	\begin{array}{cccc}
		\frac{1}{3} & 0 & 0 & 0 \\
		1 &\frac{2}{3} &-\frac{1}{3}&0
	\end{array}\right)$, then 
	$$DA-A(D\otimes I+I\otimes D)=\left(
	\begin{array}{cccc}
		\frac{a}{3}-b & -\frac{b}{3} & \frac{2 b}{3} & 0 \\
		2 a-d & \frac{2 (a+d)}{3}+b-\frac{2 d}{3} & \frac{1}{3} (-a-d)+b+\frac{d}{3} & \frac{b}{3} \\
	\end{array}
	\right).$$ Therefore $D=\begin{pmatrix}
		0 & 0\\
		c & 0
	\end{pmatrix}.$
	
	In the case $A=A_{10}(\mathbf{c})=\left(
	\begin{array}{cccc}
		0 &1 & 1 &1 \\
		\beta_1 &0 &0 &-1
	\end{array}
	\right)$ the system is
	\[\begin{array}{lc}
		2 c-b \beta_1 &=0,\\
		2 a \beta_1-d \beta_1 &=0,\\
		c+d &=0,\\
		-a+3 b+2 d &=0.\end{array}\] Hence, if $\beta_1 = 0$, then  $D=\begin{pmatrix}
		3b & b\\
		0 & 0
	\end{pmatrix}$; if $\beta_1 \neq 0,4$, then $D=0$ 
	. 
	
	In the case $A=A_{11}(\mathbf{c})=\left(
	\begin{array}{cccc}
		0 &0 & 0 &1 \\
		\beta_1 &0 &0 &0
	\end{array}
	\right)$ we have 
	$$DA-A(D\otimes I+I\otimes D)=\left(
	\begin{array}{cccc}
		-b \beta _1 & c & c & 2 d-a \\
		2 a \beta _1-\beta _1 d & b \beta _1 & b \beta _1 & -c \\
	\end{array}
	\right).$$ Hence, $D=0$, if
	$\beta_1 \neq 0;$ otherwise $D=\begin{pmatrix}
		2d & b\\
		0 & d
	\end{pmatrix}.$
	
	If $A=A_{12}(\mathbf{c})=\left(
	\begin{array}{cccc}
		0 & 1 & 1 &0  \\
		\beta_1 &0& 0 &-1
	\end{array}
	\right)$, then 
	$$DA-A(D\otimes I+I\otimes D)=\left(
	\begin{array}{cccc}
		2 c-b \beta _1 & d & d & 3 b \\
		2 a \beta _1-\beta _1 d & b \beta _1-2 c & b \beta _1-2 c & -d \\
	\end{array}
	\right).$$ Hence, $D=\begin{pmatrix}
		a & 0\\
		0 & 0
	\end{pmatrix}$, if $\beta_1=0;$  otherwise $D=0.$ 

	Here are the corresponding results in the cases of characteristic $2$ and $3.$ The proof is similar to that of the case of characteristic is not $2$ and $3.$
	
	\begin{theorem}\label{thm4} The automorphisms and derivations of the algebras, listed in Theorem 2.6, are given as follows:
		\begin{itemize}
			\item $Aut(A_{1,2}(\alpha_1, \alpha_2, \alpha_4, \beta_1)=\{I\},$\\
			$Der(A_{1,2}(\alpha_1, \alpha_2, \alpha_4, \beta_1))=\{0\},$
			\item $Aut(A_{2,2}(\alpha_1, \alpha_4,\beta_2))=\{I\}$ if $\alpha_4 \neq 0,$\\
			$Der(A_{2,2}(\alpha_1, \alpha_4, \beta_2))=\{0\},$ if $\alpha_4\neq 0$ or  $\beta_2\neq 1$, \\
			$Aut(A_{2,2}(\alpha_1,0,1))=\{\left(
			\begin{array}{cc}
				1 & 0 \\
				c & 1 \\
			\end{array}
			\right); c\in \mathbb{F}\}$,\\
			$Der(A_{2,2}(\alpha_1, 0, 1))= \{\left(
			\begin{array}{cc}
				0 & 0 \\
				c & 0 \\
			\end{array}
			\right) :\ c\in \mathbb{F}\},$
			
			\item $Aut(A_{3,2}(\alpha_1,\alpha_4,\beta_2))=
			\{\left(
			\begin{array}{cc}
				1 & 0 \\
				0 & \pm1 \\
			\end{array}
			\right)\},$ if $\alpha_4\neq 0$,\\
			$Aut(A_{3,2}(\alpha_1,0,\beta_2))=\{\left(
			\begin{array}{cc}
				1 & 0 \\
				0 & d \\
			\end{array}
			\right); d\in \mathbb{F}, d\neq 0\},$ if $\beta_2\neq 1$,\\
			$Der(A_{3,2}(\alpha_1, \alpha_4, \beta_2))=\{\left(
			\begin{array}{cc}
				0 & 0 \\
				0 & d \\
			\end{array}
			\right) :\ d\in \mathbb{F}\},$ if $\alpha_4\neq 0$ or  $\beta_2\neq 1$,\\
			$Aut(A_{3,2}(\alpha_1,0,1))=\{\left(
			\begin{array}{cc}
				1 & 0 \\
				c & d \\
			\end{array}
			\right); c,d\in \mathbb{F}, d\neq 0\},$ \\
			$ Der(A_{3,2}(\alpha_1, 0, 1))=\{\left(
			\begin{array}{cc}
				0 & 0 \\
				c & d \\
			\end{array}
			\right):\ c, d\in \mathbb{F}\},$\\  
			\item $Aut(A_{4,2}(\alpha_1,\beta_1,\beta_2))=\{I,\begin{pmatrix}1&0\\1+\beta_2&1	\end{pmatrix}\},$ \\
			$Der(A_{4,2}(\alpha_1, \beta_1, \beta_2))=\{0\},$ if $\beta_2\neq 1,$\\
			$Der(A_{4,2}(\alpha_1, \beta_1, 1))=\{\left(
			\begin{array}{cc}
				0 & 0 \\
				c & 0 \\
			\end{array}
			\right) :\ c\in \mathbb{F}\},$ 
			\item $Aut(A_{5,2}(\alpha_1,\alpha_4))=\{I\}$ if $\alpha_4 \neq 0$,\\
			$Aut(A_{5,2}(1,0))=\{ \left(
			\begin{array}{cc}
				1 & 0 \\
				c & 1 \\
			\end{array}
			\right):c \in\mathbb{F}\},$\\
			$Der(A_{5,2}(\alpha_1, \alpha_4))=\{0\},$ if $\alpha_4 \neq 0$\\
			$Der(A_{5,2}(1,0))=\begin{pmatrix}
				0 & 0\\
				c & 0
			\end{pmatrix} $
			\item $Aut(A_{6,2}(\alpha_1,\alpha_4))=\{I\}$ if $\alpha_4 \neq 0$,\\
			$Aut(A_{6,2}(\alpha_1,0))=
			\{\left(
			\begin{array}{cc}
				1 & 0 \\
				0 & d \\
			\end{array}
			\right): d\in \mathbb{F}, d\neq 0\},$ if $\alpha_1\neq 1$,\\
			$Der(A_{6,2}(\alpha_1, \alpha_4))=\{\left(
			\begin{array}{cc}
				0 & 0 \\
				0 & d \\
			\end{array}
			\right),$ if $\alpha_1 \neq 1$ or $\alpha_4 \neq 0$,\\
			$Aut(A_{6,2}(1,0))= \{\left(
			\begin{array}{cc}
				1 & 0 \\
				c & d \\
			\end{array}
			\right): c,d\in \mathbb{F}, d\neq 0\},$ \\
			$Der(A_{6,2}(1, 0))=\{\left(
			\begin{array}{cc}
				0 & 0 \\
				c & d \\
			\end{array}\right):\ c,d\in \mathbb{F}\},$
			\item $Aut(A_{7,2}(\alpha_1, \beta_1))=\{I, \begin{pmatrix}
				1 & 0 \\ 
				\alpha_1+1& 1
			\end{pmatrix}\},$\\ 
			$Der(A_{7,2}(\alpha_1, \beta_1))=\{0\},$ if $\alpha_1\neq 1, $ \\
			$Der(A_{7,2}(1, \beta_1))= \{\left( \begin{array}{cc}
				0 & 0 \\
				c & 0 \\
			\end{array}  \right):\ c\in \mathbb{F}\},$
			
			\item $Aut(A_{8,2}(\beta_1))=\\ \{I\}\cup \{\begin{pmatrix}
				1+d&1+d+d^2\\ 1&d
			\end{pmatrix}:\ \beta_1=\frac{1}{d^2+d+1}\}\cup \{\begin{pmatrix}
				d & d+d^2\\
				\frac{1+d}{d} & d 
			\end{pmatrix}:,\ \beta_1=\frac{(1+d)}{d^3}\}$, if $\beta_1 \neq 0$\\
			$Aut(A_{8,2}(0))=
			\{
			\begin{pmatrix}
				1+b & b\\
				0 & 1
			\end{pmatrix}
			:\ b\in\mathbb F
			\}.$\\
			$Der(A_{8,2}(\beta_1))=\{0\},$ if $\beta_1\neq 0$,\\
			$Der(A_{8,2}(0)=\{\begin{pmatrix}
				a & a\\
				c & c
			\end{pmatrix}:\ a,b\in\mathbb{F}\}$,\\
			
			\item		$Aut(A_{9,2}(\beta_1))=\{
			\begin{pmatrix}
				d^2 & 0\\
				0 & d
			\end{pmatrix}
			:\ d^3=1
			\}
			\cup
			\{
			\begin{pmatrix}
				0 & b\\
				\frac1b & 0
			\end{pmatrix}
			:\ \beta_1=\frac{1}{b^3}
			\},$
			if $\beta_1\neq0$,\\
			$Der(A_{9,2}(\beta_1))=\{0\}$ if $\beta_1\neq 0$,\\
			$Aut(A_{9,2}(0))
			=
			\{
			\begin{pmatrix}
				d^2 & b\\
				0 & d
			\end{pmatrix}
			:\ d\neq0,\ b,d\in\mathbb F
			\},$\\
			$Der(A_{9,2}(0))=\{\begin{pmatrix}
				0 & b\\
				0 & d
			\end{pmatrix}: \ b,d\in \mathbb{F}\},$   

			\item $Aut(A_{10,2}(\beta_1))=\{I,\left(\begin{array}{cc}
				1& 0\\
				1&1
			\end{array}\right), \begin{pmatrix}
				d & 1\\
				1+d^2 & d
			\end{pmatrix},\begin{pmatrix}
				d & 1\\
				1+d+d^2 & 1+d
			\end{pmatrix},\begin{pmatrix}
				1+d & 1\\
				1+d+d^2 & d
			\end{pmatrix},\\ \begin{pmatrix}
				1+d & 1\\
				d^2 & 1+d
			\end{pmatrix}\}$, where $\beta_1=d+d^2$,\\
			$Der(A_{10,2}(\beta_1))=\{0\},$
			\item $Aut(A_{11,2}(\beta_1))=\{\left(\begin{array}{cc}
				a & 0\\
				c & 1
			\end{array}\right):\ 0\neq a,c\in\mathbb{F},  c^2=(1+a^2)\beta_1\}$, \\ 
			$Der(A_{11,2}(\beta_1))=\{\left(
			\begin{array}{cc}
				a & 0 \\
				c & 0 \\
			\end{array} \right):\ a,\ c \in \mathbb{F}\}.$
			
		\end{itemize}
	\end{theorem}
	
		\begin{theorem}\label{thm4} The automorphisms and derivations of the algebras, listed in Theorem 2.7, are given as follows:
		\begin{itemize}
			\item $Aut(A_{1,3}(\alpha_1,\alpha_2,\alpha_4,\beta_1))=\{I\},$\\
			$Der(A_{1,3}(\alpha_1, \alpha_2, \alpha_4, \beta_1))=\{0\},$
			\item $Aut(A_{2,3}(\alpha_1, \alpha_4, \beta_2))=\{I\}$,\\
			$Der(A_{2,3}(\alpha_1, \alpha_4, \beta_2))= \{0\},$
			\item  
			$Aut(A_{3,3}(\alpha_1,\alpha_4, \beta_2))=\{\left(\begin{array}{cc}
				1 & 0\\
				0 & \pm 1
			\end{array}\right)\}$, if $\alpha_4\neq 0$, \\
			$Der(A_{3,3}(\alpha_1, \alpha_4, \beta_2))=\{0\},$ if $\alpha_4 \neq 0$,\\
			$Aut(A_{3,3}(\alpha_1,0, \beta_2))=\{\left(\begin{array}{cc}
				1& 0\\
				0&d
			\end{array}\right):   d\in\mathbb{F}^*\}$, if $\beta_2 \neq 2 \alpha_1-1,$\\
			$ Der(A_{3,3}(\alpha_1, 0, \beta_2))=\{\left(
			\begin{array}{cc}
				0 & 0 \\
				0 & d \\
			\end{array}
			\right)
			:\ d\in \mathbb{F}\},$ if $\beta_2 \neq 2 \alpha_1-1$,\\
			$Aut(A_{3,3}(\alpha_1,0, 2 \alpha_1-1))=\{\left(\begin{array}{cc}
				1& 0\\
				c&d
			\end{array}\right): \ c,d\in\mathbb{F}, d\neq 0\}$, \\
			$Der(A_{3,3}(\alpha_1, 0, 2 \alpha_1-1))=\{\left(
			\begin{array}{cc}
				0 & 0 \\
				c & d \\
			\end{array}
			\right)
			:\ c,\ d\in \mathbb{F}\},$ 
			\item $Aut(A_{4,3}(\beta_1,\beta_2))=\{I\}$,\\
			$Der(A_{4,3}(\beta_1, \beta_2))=\{0\},$
			
			\item $Aut(A_{5,3}(\alpha_1))=\{\left(
			\begin{array}{cc}
				1 & 0 \\
				c & 1 \\
			\end{array}
			\right)\}$,\\
			$Der(A_{5,3}(\alpha_1))=\{\left(
			\begin{array}{cc}
				0 & 0 \\
				c & 0 \\
			\end{array} \right):\ c\in \mathbb{F}\},$
			\item $Aut(A_{6,3}(\alpha_1,\alpha_4))=\{I\}$,\\
			$Der(A_{6,3}(\alpha_1,\alpha_4))=\{0\}$,
			\item  $Aut(A_{7,3}(\alpha_1, \alpha_4))=\{\left(
			\begin{array}{cc}
				1 & 0 \\
				0 & \pm 1 \\
			\end{array}
			\right)\}$,  if $\alpha_4 \neq 0$,\\
			$Der(A_{7,3}(\alpha_1, \alpha_4))=\{0\},$ if $\alpha_4 \neq 0,$\\
			$Aut(A_{7,3}(\alpha_1, 0))=\{\left(
			\begin{array}{cc}
				1 & 0 \\
				0 & d \\
			\end{array}
			\right):\ 0\neq d\in\mathbb{F}\}$, \\
			$Der(A_{7,3}(\alpha_1,0))=\{\left(
			\begin{array}{cc}
				0 & 0 \\
				0 & d \\
			\end{array} \right):\ d\in \mathbb{F}\},$ 
			
			\item $Aut_{8,3}(\beta_1)=\{I\}$,\\
			$Der(A_{8,3}(\beta_1))=\{0\},$
			\item $Aut(A_{9,3}(\beta_1))
			=
			\{I\}
			\cup\{\begin{pmatrix}
				-d & \frac{-d(1+d)}{1+2d} \\ \frac{2d^2-d-1}{d} & d
			\end{pmatrix}: \beta_1=\frac{d^3-1}{d^3}, d \neq 0,1 \}$ for $\beta_1\neq 0,1,2,$\\
			$Aut(A_{9,3}(0))=\{\begin{pmatrix}
				1 & b \\ 0 & 1
			\end{pmatrix}, b \in \mathbb{F} \},$\\
			$Aut(A_{9,3}(1))=\{\begin{pmatrix}
				\dfrac{1+d}{2} & \dfrac{1-d}{2}\\
				1-d & d
			\end{pmatrix}:d \in \mathbb{F}\},$\\
			$Aut(A_{9,3}(2))=\{I\cup\begin{pmatrix}
				1 & 0 \\ 1 & 2
			\end{pmatrix} \},$\\
			$Der(A_{9,3}(\beta_1)=\begin{pmatrix}
				2 b \beta_1 & b\\
				2 b \beta_1 & b \beta_1
			\end{pmatrix}:\ b\in\mathbb{F}\}$ for $\beta_1\in\mathbb{F},$\\
			\item 	$Aut(A_{10,3}(\beta_1))=\{\left(
			\begin{array}{cc}
				a & 0 \\
				0 & a^2 \\
			\end{array}
			\right):\ a^3=1\}\cup \{\begin{pmatrix}0&b\\  
				1/b&0\end{pmatrix}:\ b^{-3}=\beta_1\},$ if $\beta_1\neq 0$,\\
			$Der(A_{10,3}(\beta_1))=\{\begin{pmatrix}
				a & 0\\
				0 & -a
			\end{pmatrix}: \ a\in \mathbb{F}\}$ if $\beta_1\neq 0$,\\
			$Aut(A_{10,3}(0))=\{\left(
			\begin{array}{cc}
				d^2 & b \\
				0 & d \\
			\end{array}
			\right):\ b,d\in \mathbb{F},\ d\neq 0\},$\\
			$Der(A_{10,3}(0))=\{\begin{pmatrix}
				a & b\\
				0 & -a
			\end{pmatrix}: \ a,b\in \mathbb{F}\},$   
			
			\item 
			$Aut(A_{11,3}(\beta_1))= \{\begin{pmatrix}
				\pm1&0\\ 0&1
			\end{pmatrix}\},$ if $\beta_1\neq 0$,\\
			$Der(A_{11,3}(\beta_1))=\{\left(
			\begin{array}{cc}
				0 & \frac{2 c}{\beta_1} \\
				c & 0 \\
			\end{array} \right):\ c \in \mathbb{F}\},$ if $\beta_1\neq 0,$\\
			$Aut(A_{11,3}(0))= \{\begin{pmatrix}
				a&b\\ 0&1
			\end{pmatrix}\ : 0\neq a\in \mathbb{F}\},$ if $\beta_1=0$,\\
			$Der(A_{11,3}(0))=\{\left(
			\begin{array}{cc}
				a & b \\
				0 & 0 \\
			\end{array} \right):\ a, \ b \in \mathbb{F}\}.$

		\end{itemize}
	\end{theorem}

	\end{document}